\documentclass{amsart}

\title[Non-Autonomous Basins]{Non-Autonomous Basins of Attraction with $4$-Dimensional Boundaries}
\author{Han Peters, Erlend Forn\ae ss Wold}
\date{October 3, 2004}

\usepackage{amscd}

\newtheorem{theorem}{Theorem}

\newtheorem{corollary}{Corollary}

\theoremstyle{definition}

\newtheorem{conjecture}{Conjecture}

\theoremstyle{remark}
\newtheorem{remark}{Remark}

\newcommand{\NN}{\mathbb{N}}
\newcommand{\RR}{\mathbb{R}}
\newcommand{\CC}{\mathbb{C}}

\newcommand{\aut}{\mathrm{Aut}(\CC^k)}
\newcommand{\boxd}{\overline{\dim}_B}

\def\B{{\bf B}}

\def\e{{\epsilon}}
\def\d{{\delta}}

\def\r{{\rho}}

\begin{document}

\begin{abstract}
We study whether the basin of attraction of a sequence of
automorphisms of $\CC^k$ is biholomorphic to $\CC^k$. In
particular we show that given any sequence of automorphisms with
the same attracting fixed point, the basin is biholomorphic to
$\CC^k$ if the maps are repeated often enough. We also construct
Fatou-Bieberbach domains whose boundaries are $4$-dimensional.
\end{abstract}

\maketitle

\section{Introduction}

In the 1920's Fatou and Bieberbach proved the existence of proper
subdomains of $\CC^2$ that are biholomorphically equivalent to
$\CC^2$, later known as Fatou-Bieberbach domains. Their examples
were the basins of attraction of some automorphisms of $\CC^2$
which have more than one fixed point. In fact, the basin of
attraction of an attracting fixed point of an automorphism of
$\CC^k$ is always biholomorphic to $\CC^k$, which follows from the
work of Sternberg in 1957 \cite{sb} and was proved independently
by Rosay and Rudin in 1988 \cite{rr}.

More recently, Stens\o nes showed \cite{st} that there exist
Fatou-Bieberbach domains whose boundaries are smooth, so in
particular they have Hausdorff dimension $3$. Wolf \cite{wo}
showed that for any $h \in (3,4)$, there exists a Fatou-Bieberbach
domain whose boundary has Hausdorff dimension $h$. The results of
Stens\o nes and Wolf leave open two questions about the Hausdorff
dimension of the boundary a Fatou Bieberbach domain:

1) Is it possible for the dimension of the boundary to be less
than $3$?

2) Is it possible that the dimension of the boundary is exactly
$4$?

The Fatou-Bieberbach domains that were constructed in the paper by
Wolf are basins of attraction of some particular polynomial
automorphisms of $\CC^2$. We will show that one can use a sequence
of these same maps to get a Fatou-Bieberbach domain whose boundary
has upper-box dimension $4$. Since the Hausdorff dimension can in
general be larger than the upper-box dimension, this does not
answer question (2) above.

We will show however that one can use a more direct approach to
obtain the following theorem:

\begin{theorem} \label{hausdorff} There exists a
Fatou-Bieberbach domain $\Omega$ so that $\partial \Omega =
\partial \overline{\Omega}$, and such that the $4$-dimensional
Hausdorff measure of $\partial \Omega$ near any point is non-zero.
\end{theorem}

The Fatou-Bieberbach domain in the above theorem will also be a
basin of attraction of a sequence of mappings. The boundary with
upper-box dimension $4$ is constructed as the limit of lower
dimensional sets, but to prove Theorem \ref{hausdorff} we will
start with sets of Hausdorff dimension $4$ and then construct our
maps in such a way that these sets lie in the boundary of the
basin of attraction.

We will need two theorems that tell us that the basins we
construct are really biholomorphic to $\CC^k$. For a given
sequence of automorphisms of $\CC^k$ which all have the same
attracting fixed point, say $0$, define the basin of attraction by

$$
\{z \in \CC^k \mid f_n \circ \ldots \circ f_1 (z) \rightarrow 0\}.
$$

One could ask whether such a basin is always biholomorphic to
$\CC^k$. In general the answer to this question is no, which is
easy to see. However, it is possible to prove that the basin is
always biholomorphic to $\CC^k$ if one puts some restrictions to
the sequence of maps. For instance one can prove the following
theorem:

\begin{theorem}\label{square}
Let $f_1, f_2, \ldots$ be a sequence of automorphisms of $\CC^k$.
Suppose that there exist $0 < a < b < 1$ with $b^2 < a$ such that
every $f_n$ satisfies
$$
a\|z\| \le \|f_n(z) \| \le b \|z\|,
$$
for every $z$ in the unit ball. Then the basin of attraction of
the sequence is biholomorphic to $\CC^k$.
\end{theorem}
The proof of this theorem can be found in \cite{wd}. We will use
this theorem to prove Theorem \ref{hausdorff}.

Though the above theorem is very useful for the construction of
Fatou-Bieberbach domains, it is somewhat unsatisfying because all
the maps must satisfy very strict conditions on the speed at which
orbits converge to $0$. We will prove the following generalization
of the main theorem of the appendix in \cite{rr}:

\begin{theorem} \label{repeat}
Let $f_1, f_2, \ldots$ be a sequence of automorphisms of $\CC^k$
which all have the same attracting fixed point. Then we can find
large enough constants $n_1, n_2, \ldots $ such that the basin of
attraction of the sequence $f_1^{n_1}, f_2^{n_2}, \ldots$ is
biholomorphic to $\CC^k$.
\end{theorem}

We show how one needs to alter the definition of the basin of
attraction of a sequence for this theorem to hold. We will then
show how one can use Theorem \ref{repeat} to obtain a
Fatou-Bieberbach domain whose boundary has upper-box dimension $4$
using the maps used in \cite{wo}.

Basins of attraction of sequences of mappings have an interesting
connection to complex dynamical systems. The following conjecture
was posed by Bedford \cite{be}:

\begin{conjecture}\label{bedford}
Let $F$ be an automorphism of a complex manifold, which is
hyperbolic on a compact set $K$. Then for every $p \in K$ the
stable manifold is biholomorphically equivalent to complex
Euclidean space.
\end{conjecture}

So far the best answer to this conjecture has been given by
Jonsson and Varolin \cite{jv}, who proved that with respect to any
invariant probability measure, almost every stable manifold is
biholomorphic to $\CC^k$. It follows from the work of Forn\ae ss
and Stens\o nes \cite{fst} that the above conjecture can be
answered positively by proving the following conjecture about
non-autonomous basins of attraction:

\begin{conjecture}\label{compact}
Let $f_1, f_2, \ldots$ be a sequence of automorphisms of $\CC^k$
and assume that there exist $0 <a < b < 1$ such that for every $n$
and every $z$ in the unit ball the following holds:
$$
a\|z\| \le \|f_n(z)\| \le b \|z\|.
$$
Then the basin of attraction of $0$ for this sequence is
biholomorphically equivalent to $\CC^k$.
\end{conjecture}

By using the biholomorphisms that map the local stable manifolds
onto the local stable tangent bundle, one can translate the
setting of Conjecture \ref{bedford} into a sequence of
biholomorphic mappings from the unit ball into the unit ball,
which satisfy the conditions of Conjecture \ref{compact}. In
\cite{fst} it is shown how one can define the basin of attraction
of such a sequence of biholomorphic mappings in terms of the
{\emph{tail space}}. It is easy to see that this basin is
biholomorphic to the stable manifold. It is also shown that this
basin of attraction is biholomorphic to the basin of attraction of
sequence of global automorphisms which satisfy the conditions in
Conjecture \ref{compact}. To extend the local biholomorphic
mappings to the global automorphisms the following theorem, which
was proved by Forstneric \cite{fc} and independently by Weickert
\cite{we}, was used:

\begin{theorem}\label{automorphism}
Let $P = (P_1, \ldots , P_k)$, $k \ge 2$ be a holomorphic
polynomial mapping of $\CC^k$ to itself, with $P^\prime(0)$
invertible. Let $d \ge max_i(deg(P_i))$. Then there exists $\phi
\in \mathrm{Aut}(\CC^k)$ such that the $d$-jet of $\phi$  at $0$
equals $P$.
\end{theorem}

We will use Theorem \ref{automorphism} in the proof of theorem
\ref{repeat}.

In the second section we will introduce the notation that we will
use throughout the article. In the third section we will prove
Theorem \ref{repeat}. In section four we will look at two
interesting examples of basins of attractions that are not
biholomorphic to complex Euclidean space. These examples
demonstrate that Theorem \ref{repeat} does not hold when we loosen
the conditions. In the fifth section we will construct
Fatou-Bieberbach domains whose boundaries have upper box dimension
$4$. In the sixth section we will prove Theorem \ref{hausdorff}.
In the last section we will give an elementary proof of the fact
that the Hausdorff dimension of the boundary of a Fatou-Bieberbach
domain that is Runge is always at least $3$.

\verb"Acknowledgement" We would like to thank Christian Wolf for
suggesting that a Fatou-Bieberbach domain with dimension-4
boundary can be constructed as the basin of attraction of a
sequence of H\'{e}non maps. It was his suggestion that got us
started on the subject.

\section{Preliminaries}

From now on we fix an integer $k\ge2$, and we will write $\aut$
for the automorphisms of $\CC^k$ that fix the origin. Throughout
the whole paper we will write products and powers for the
composition of maps, e.g. $fg$ and $f^2$ for $f\circ g$ and $f
\circ f$.

Throughout the paper we will work with sequences of maps $f_1,
f_2, \ldots \in \aut$ and study orbits $z , f_1(z), f_2 f_1(z),
\ldots$. For $f(n)(z) = f_n  \cdots f_1(z)$ we will often write
$z_n$.

When also a sequence of integers $n_1, n_2, \ldots$ is given and
the maps $f_j$ are iterated $n_j$ times, then we will write $F_j =
f_j^{n_j}$. We will also use the notation $N_j = n_j + \cdots +
n_1$, and we will write
$$
F(j) = F_j \cdots F_1.
$$
We think of $F(j)$ as the map that takes us from stage $0$ to
stage $j$. For some $m \ge n$, we will write $F(n,m)$ for the map
that takes us from stage $n$ to stage $m$, i.e.
$$
F(n,m) = F(m) F(n)^{-1}.
$$
For $n=N_j +m$ (where $m < n_{j+1}$) we will also write
$$
f(n) = f_{j+1}^m F(j) \; \mathrm{and} \; f(p, q)= f(q)f(p)^{-1}.
$$
Notice that $f(n)$ is the composition of $n$ maps $f_j$, while
$F(n)$ is the composition of $n$ maps $F_j$.

For $z \in \CC^k$ we will write $\|z\|$ for the Euclidean norm of
$z$. We will write $\B$ for the unit ball in $\CC^k$, and for $z
\in \CC^k$ and $r>0$ we write $B(z, r)$ for the ball of radius $r$
centered at $z$. We will also write $B(r)$ for $B(0,r)$.

Let $f_1, f_2, \ldots \in \aut$ and assume that every $f_n$ is
attracting at $0$. The usual definition for the basin of
attraction of this sequence is
$$
\Omega^\star = \{ z \in \CC^k \mid f(j)(z) \rightarrow 0\}.
$$
This definition does not work well for our purposes. Suppose that
there exist radii $r_1, r_2,\ldots$, and some $\rho < 1$, such
that $\|f_j(z)\| < \rho \|z\|$ for any $z \in \Delta(r_j)$, and
such that $f_j(B(r_j)) \subset B(r_{j+1})$. We define the basin of
attracting of $0$ for the sequence $f_1, f_2, \ldots$ as
$$
\Omega = \bigcup f(j)^{-1} B(r_{j+1}).
$$
Notice that the definition of the basin depends on the choices of
the $r_j's$. It is easy to see that if all the radii $r_j$ are
equal then $\Omega$ is equal to $\Omega^\star$. However, if the
radii $r_j$ shrink then there might be orbits that converge to $0$
but that never reach the balls $B(r_j)$ at the appropriate stage.
We will see later that this can indeed happen and that in such
cases $\Omega$ might be biholomorphic to $\CC^k$ while
$\Omega^\star$ is not.  If we want to put emphasis on which
sequence we are referring to, we will denote the basin of
attraction by $\Omega_{\{f_j\}}$.

One might want to compare our definition for the basin of
attraction with the definition of the basin of attraction for a
semi-attracting fixed point, see for instance \cite{ue}.

\section{Basins Biholomorphic to $\CC^k$}

In this section we will be working with invertible maps $T_j, g_j$
and $G_j$, where $G_j = g_j^{n_j}$ and we will use the notation
$$
G(j) = T_j^{-1}G_j T_j \cdots T_1^{-1} G_1 T_1, \; \mathrm{and} \;
G(n,m) = G(m) G(n)^{-1},
$$
and for $ n = N_j + m$ we write
$$
g(n) = T_{j+1}^{-1} g_j^m T_{j+1} G(j) \; \mathrm{and} \: g(p,q) =
g(q)g(p)^{-1}.
$$
Usually we will be working with inverse mappings for the $G$'s, so
then we can think of $G(n,m)^{-1}$ as the map that takes us back
from stage $m$ to stage $n$. The definition for $G(n)$ may not be
entirely consistent with the definition for $F(n)$, but this will
not cause any problems.

We will now prove Theorem \ref{repeat}:\\
\\
{\bf{Theorem \ref{repeat}}} {\emph {Let $f_1, f_2, \ldots$ be a
sequence of automorphisms of $\CC^k$ which all have the same
attracting fixed point. Then we can find large enough constants
$n_1, n_2, \ldots $ such that the basin of attraction of the
sequence $f_1^{n_1}, f_2^{n_2}, \ldots$ is biholomorphic to
$\CC^k$.}}
\\

By large enough we mean that we may replace any constant $n_j$ by
a larger number, but then we may also have to enlarge all
subsequent constants $n_{j+1}, n_{j+2}, \ldots$.

\begin{remark} \label{ready}
For the basin to be well-defined we choose any sequence of radii
$r_j$ such that $B(r_j)$ is contained in the basin of $f_j$. We
may replace $f_j$ by a high iterate to get that $f_j(B(r_j))
\subset \subset B(r_{j+1})$ and $\|f_j(z)\| < \frac{1}{2} \|z\|$
for $z \in B(r_j)$. It is easy to see that the theorem does not
really depend on the choice of $r_j's$, since for two different
sequences one can increase the $n_j's$ such that the basins are
identical. We may decrease the values of the $r_j$'s in the proof,
but we can always increase the values of the $n_j$'s later such
that the definition of the basin of attraction with the smaller
$r_j$'s is equivalent to the definition of the basin of attraction
with the original $r_j$'s.
\end{remark}

\begin{proof}
For simplicity we assume that $0$ is the attracting fixed point.

It follows from Lemma 3 in the appendix of \cite{rr} that for
every $j \in \NN$ we can find:

(1) A polynomial automorphisms $g_j$ of $\CC^k$ which is linearly
conjugate to a lower triangular automorphism, with $g_j(0)=0$ and
$g_j^\prime(0) = f_j^\prime(0)$, and

(2) For $m_j$ as large as we desire a polynomial map $\phi_j:\CC^k
\rightarrow \CC^k$, with $\phi_j(0)=0, \phi_j^\prime(0) = I$, and
for which the following equation holds
$$
g_j^{-1} \phi_j  f_j - \phi_j = O(\|z\|^{m_j}).
$$
It follows from Theorem \ref{automorphism} that for every $j \in
\NN$ we can find an automorphisms $T_j$ with $T_j - \phi_j =
O(\|z\|^{m_j})$. Therefore we get the equation:
$$
T_j^{-1} g_j^{-1}  T_j  f_j - I = O(\|z\|^{m_j}).
$$
Recall from \cite{rr} that there exist constants $\gamma_j$ such
that
$$
\|g_j^{-n}(w) - g_j^{-n}(w^\prime)\| \le \gamma_j^n \|w -
w^\prime\|,
$$
for all $w, w^\prime$ in $\B$ and all natural numbers $n$.  Also
recall that the $g_j$'s are attracting at $0$, and that the basins
of the maps $g_j$ are all of $\CC^k$.

We noted in Remark \ref{ready} that we may assume that $\|f_j(z)\|
\le \frac{1}{2} \|z\|$ for all $j \in \NN$ and all $z \in B(r_j)$.
We now choose the $m_j$'s so large that
$$
\frac{\gamma_j}{2^{m_j}} <1,
$$
and such that $m_{j+1} > m_j$ for every $j \in \NN$. Then it
follows from equation (6) of the proof of the theorem in the
appendix of \cite{rr} that there exists constants $C_j$ such that
for every $j$ and every $z \in B(r_j)$ the following holds
$$
\|T_j^{-1}g_j^{-n-1}T_jf_j^{n+1}(z)
-T_j^{-1}g_j^{-n}T_jf_j^n (z)\| \le C_j
(\frac{\gamma_j}{2^{m_j}})^n \|z\|^{m_j}.
$$
Since $(\frac{\gamma_j}{2^{m_j}})^n$ is summable and the other
terms don't depend on $n$, we can decrease the radii $r_j$ if
necessary such that for $z \in B(r_j)$ we have
\begin{eqnarray}\label{sum}
\sum_{n \ge 0} \|T_j^{-1}g_j^{-n-1}T_jf_j^{n+1}(z)
-T_j^{-1}g_j^{-n}T_jf_j^n (z)\| \le \frac{1}{2} \|z\|^{m_j - 1}.
\end{eqnarray}
It follows that for $z \in B(r_j)$ and any $m$ larger than $n$ we
have
\begin{eqnarray} \label{convergence}
\|T_j^{-1}g_j^{-m}T_j f_j^{m}(z) -T_j^{-1} g_j^{-n} T_j f_j^n
(z)\| \le \frac{1}{2} \|z\|^{m_j - 1}.
\end{eqnarray}
We now inductively construct large radii $R_j>1$ and large
integers $n_j$. Suppose that we have constructed $R_1$ through
$R_{j-1}$ and $n_1$ through $n_{j-1}$. At this stage we have
already fixed the automorphism $G(j-1)$, therefore we can fix
$R_j$ such that
\begin{align} \label{large}
G(j-1)(B(2^j))\subset B(R_j).
\end{align}

Then choose $n_j$ at least so large that
\begin{eqnarray} \label{shrinking}
f_j^{n_j}(B(r_j)) \subset B(r_{j+1}),
\end{eqnarray}
and
\begin{eqnarray} \label{growing}
T_j(B(R_j)) \subset g_j^{-n_j}T_j(B(\frac{r_{j+1}}{2})).
\end{eqnarray}

It follows from \eqref{sum}, \eqref{convergence} and
\eqref{shrinking} and the fact that we choose $m_{j+1}$ strictly
larger than $m_j$ that we can increase $n_j$ if necessary and
choose all subsequent $n_{j+1}, n_{j+2}, \ldots$ large enough such
that the following inequality will hold throughout the
construction
\begin{eqnarray} \label{final}
\|g(N_{j-1}, m)^{-1} f(N_{j-1}, m)(z) - z \| \le \|z\|^{m_j - 1},
\end{eqnarray}
for every for every $m$ and every $z \in B(r_j)$.

We will now show that with these choices the basin $\Omega$ is
biholomorphic to $\CC^k$.

Let $K_j = F(j)^{-1}(B(r_{j+1}))$. Clearly $F(j)(K_j) =
B(r_{j+1})$, and thus we have that
$$
\|g(N_j, m)^{-1}f(m)(z) - F(j)(z)\| \le \|F(j)(z)\|^{m_{j+1}-1},
$$
for every $z \in K_j$. It follows that $H_m = g(m)^{-1} f(m)$ is
bounded on every $K_j$. Therefore we can find a subsequence of
$\{H_m\}$ that converges uniformly on $K_1$, say to a map $H$. We
have that $K_1 \subset K_2 \subset \ldots$ and thus we can apply
the same argument to every $K_j$, and by taking a smaller and
smaller subsequences we see that $H$ extends to $\bigcup K_j =
\Omega$. We claim that $H$ maps $\Omega$ biholomorphically onto
$\CC^k$.

It is a well known fact that the limit map of a sequence of
biholomorphic maps is either degenerate everywhere or it is one to
one. $H$ is one to one since $H_m^\prime(0) = I$ for all $m$.

We still need to show that $H$ is surjective. Notice that it
follows from \eqref{convergence} and \eqref{shrinking} that
$$
g(N_j, m)^{-1} f(m) (K_j) \supset B(\frac{r_{j+1}}{2}),
$$
for all $m$. Therefore, it follows from \eqref{large} and
\eqref{growing} that $H_n(K_j)$ contains the ball of radius $2^j$
for all $n \ge N_j$, and thus $B(2^j) \subset H(K_j)$. Therefore
$H$ maps $\Omega$ biholomorphically onto $\CC^k$.
\end{proof}

We will later use this theorem to construct a Fatou-Bieberbach
domain whose boundary has upper box dimension $4$.

\section{Counterexamples}

In this section we present some examples that show that Theorem
\ref{repeat} is as good as it can get. We first give an example by
Forn\ae ss that shows that the result does not hold if we don't
allow our maps to repeat. A second example will show that the
theorem does not hold if we choose a more straightforward
definition for $\Omega$.

Recall the following result from \cite{fo}
\begin{theorem} \label{fornaess}
Let $a_n$ be complex numbers such that $0 < |a_n| < 1$ and such
that $|a_{n+1}| < |a_n|^t$ for some fixed $t>2$. Define the
polynomial automorphisms $f_j(z,w) = (z^2 + a_j w, a_j z)$ and
define
$$
\Omega = \{ z \in \CC^2 \mid f(n) (z) \rightarrow 0\}.
$$
Then $\Omega$ is not biholomorphic to $\CC^2$
\end{theorem}

It is shown in \cite{fo} that the basins from theorem
\ref{fornaess} are increasing unions of balls, but also not
biholomorphic to the ball or the cylinder. Therefore they are
examples of so-called short $\CC^2$'s. The definition of $\Omega$
in theorem \ref{fornaess} is equivalent to ours since all the maps
are contracting on one fixed ball around the origin, thus proving
that our theorem does not hold if we don't allow the maps to
repeat.

Let us recall some simple facts about H\'{e}non mappings in
$\CC^2$, that can for instance be found in \cite{bs}. Let $F \in
\mathrm{Aut}(\CC^2)$ be of the form
$$
F(z,w) = (aw + P(z), z),
$$
with non-zero $a \in \CC$ and $P$ a polynomial of degree at least
$2$. For $R>0$ one defines:
\begin{align*}
V^+ = \{(z, w) \in \CC^2 \mid |z| \ge |w|, |z| \ge R \},\\
V^- = \{(z, w) \in \CC^2 \mid |w| \ge |z|, |w| \ge R \},
\end{align*}
and
$$
D = \{(z,w) \in \CC^2 \mid |z|, |w| \le R\}.
$$
Then one can choose $R$ large enough such that the following
properties hold:
\begin{align}
F(V^+) \subset V^+, \\
F^{-1}(V^-) \subset V^-,
\end{align}
and the orbit of some $z \in \CC^2$ converges to infinity in
positive time (or negative time) if and only if for some $n \in
\NN$ one has that $F^n(z) \in V^+$ (resp. $F^{-n}(z) \in V^-$).
The sets $D, V^+ \ \mathrm{and} \ V^-$ together are called a
\emph{filtration} for $F$. The maps that we use in this paper are
usually of the form $F(z,w) = (aw +z^2, az)$, but this map is just
a linear conjugation of a H\'{e}non map, so there also exist
filtration for such a mapping.

 For the following theorem, define the maps $f_j = (z^2 + a_j
w, a_j z)$ as for Theorem \ref{fornaess}, but let the constants
$a_j$ be real numbers in $(0,1)$ that converge to $1$. Recall that
$\Omega^\star$ is defined as $\{z \in \CC^2 \mid f(n)(z)
\rightarrow 0\}$.

\begin{theorem}\label{slowing}
$\Omega^\star$ is not open, and is therefore in particular not
biholomorphic to $\CC^2$.
\end{theorem}

\begin{proof}
Let $\Omega_\infty$ be the set of all points in $\CC^2$ which
escape to infinity. Since the $a_n$ are all smaller than $1$,
there exist a uniform filtration for all the mappings $f_j$, in
fact $R = 3$ will suffice. Therefore we have that if $z \in
\Omega_\infty$, then there is some $z_n = (x_n, y_n)$ for which
$|x_n|
> 3$ and $|x_n| > |y_n|$. Also, if $z=(x,y)$ is such that $|x| >
3$ and $|x| > |y|$, then we have that $\|f_n(z)\| > \|z\| +1$.
Therefore it follows that
\begin{eqnarray*}
\Omega_\infty = \bigcup_{n \ge 0} f(n)^{-1}\left(\{z=(x,y) \in
\CC^2 \mid |x|>3, |x| > |y|\} \right).
\end{eqnarray*}
We see that $\Omega_\infty$ is a union of open sets, and thus
open. It is not clear whether the complement of $\Omega_\infty$ is
equal to $\Omega^\star$ or not. However, we shall see that
restricted to the first quadrant of $\RR^2$ the set
$\Omega_\infty$ is exactly the complement of $\Omega^\star$.
Therefore, $\Omega \cap \RR^2_+$ is closed and non-empty, so
$\Omega^\star$ is not open in $\CC^2$.

Let $z = (x,y)$ be a point in $\RR^2_+:=\{(x,y) \mid x,y \ge 0\}$.
We will write $\|z\|_1$ for $x + y$. Assume that $z_n$ does not
converge to $0$. Then there exists an $\epsilon> 0$ such that
$\|z_n\|_1 > \epsilon$ for arbitrary large $n$. We take $j$ so
large that for all $n \ge j$ we have
$$
(a_{n})^2(1+\frac{1}{4}\epsilon^2) > 1 + \frac{1}{5}\epsilon^2,
$$
and we choose $n \ge j$ such that $\|z\| > \epsilon$. We have
$$
f_{n+2} f_{n+1}(z_n) = ((x_n^2 + a_{n+1}y_n)^2 + a_{n+2}a_{n+1}
x_n, a_{n+2} x_n^2 + a_{n+2}a_{n+1} y_n ),
$$
and therefore
\begin{align*}
\|f_{n+2} f_{n+1}(z_n)\|_1 &\ge a_{n+1}^2 y_n^2 + a_{n+2}x_n^2 + a_{n+2}a_{n+1}(x_n+y_n)\\
&> \min(a_{n+1}, a_{n+2})^2(x_n+y_n+x_n^2+y_n^2)\\
&\ge \min(a_{n+1}, a_{n+2})^2(1+\frac{1}{4}\epsilon^2) \|z\|_1 >
(1+\frac{1}{5}\epsilon^2)\|z\|_1
\end{align*}
Therefore $\|z_n \|_1$ goes to infinity and $z \in \Omega_\infty$.
Hence we have
\begin{eqnarray*}
\RR^2_+ = (\Omega^\star \cap \RR^2_+) \cup (\Omega_\infty \cap
\RR^2_+).
\end{eqnarray*}
We have that $\Omega_\infty \cap \RR^2_+$ is a nonempty relatively
open subset of $\RR^2_+$, and therefore $(\Omega^\star \cap
\RR^2_+)$ is a proper relatively closed subset of $\RR^2_+$. Since
$0\in \Omega$, it follows that $\Omega^\star$ can't be open.
\end{proof}

This shows that Theorem \ref{repeat} does not hold if we use the
more standard definition $\Omega^\star$ for the basin of
attraction, even if we allow the maps to repeat as often as we
want.

\begin{remark}
If one takes the maps $f_n(z,w) = (a_n w + z^2, a_n w)$ and
requires that there exist some constants $0<a<b<1$ such that for
every $n \in \NN$ one has that $a \le a_n \le b$, then the
conditions in Conjecture \ref{compact} are satisfied. In this case
it is easy to see that the basin of attraction is biholomorphic to
$\CC^2$. One just takes $A_n = f_n^{\prime}(0)$ and defines $H_n =
A(n)^{-1} f(n)$, where $A(n) = A_n \ldots A_1$. Since for $z$
close to $0$ one has that $\|f_n(z)\| \le a_N^{1.1}\|z\|$, one
gets that $\{H_n(z)\}$ is a Cauchy sequence for $z \in \Omega$,
and the maps $H_n$ converge uniformly on compact subsets of
$\Omega$. The limit mapping is a biholomorphic mapping from
$\Omega$ onto $\CC^k$.
\end{remark}

\section{Upper Box Dimension $4$}

We will start this section by recalling the definition of upper
box dimension. Let $K$ be a compact subset of some metric space.
For $\epsilon >0$ we write $\mathcal{B}_\epsilon$ for the set of
all coverings $\{B_i\}$ of $K$ with balls of radius $\epsilon$.
For $h
> 0$ we define
$$
\gamma_h^\epsilon(K) = \epsilon^h \inf_{\mathcal{B}_\epsilon} \#
\{B_i\}, \; \mathrm{and}
$$
$$
\mu_h(K) = \limsup_{\epsilon \rightarrow 0}
\gamma_h^{\epsilon}(K).
$$
$\mu_h(K)$ is called the $h$-upper box content (or upper Minkowski
content) of $K$.

For every set $K$ there is exactly one value $h$ such that
$\mu_{h^\prime}(K) = 0$ for all $h^\prime > h$ and
$\mu_{h^\prime}(K) = \infty$ for all $h^\prime < h$. This $h$ is
called the upper box dimension of $K$, and is written by
$\boxd(K)$. The upper box dimension of a set may be larger than
its Hausdorff dimension, but often the two concepts are equal.
Upper box dimension is often used instead of Hausdorff dimension
because it is better suited for computer approximations.

For two compact subsets $A,B \subset \CC^k$ we use the usual
definition for the Hausdorff distance between $A$ and $B$, namely
$$
d_H(A,B) = \max\{d(A,B) , d(B, A)\},
$$
where
$$
d(A,B) = \sup_{x\in A}\inf_{y\in B}d(x,y).
$$

Before we continue to construct the Fatou-Bieberbach domain whose
boundary has upper-box dimension $4$ we first need to recall some
more facts about H\'{e}non maps, which can be found in \cite{bs}.
We have already defined the filtration $D$, $V^+$, $V^-$ for a
H\'{e}non map $F$. We now define the sets $K^+$ and $K^-$ as
$$
K^{\pm} = \{z \in \CC^k \mid \{F^{\pm n}(z)\} \mathrm{\> is \>
bounded } \}.
$$
It is easy to see that $z \in K^+$ if and only if $F^{n}(z) \in D$
for all $n$ large enough.

We define the forward and backward Julia sets $J^+$ and $J^-$ as
$J^\pm = \partial K^\pm$, and the Julia set $J = J^+ \cap J^-$.

The dynamics of conjugated H\'{e}non maps of the form $f(z,w) =
(z^2 + c + a w, a z)$, where $c$ in the main cardioid of the
Mandelbrot set and the constant $a>0$ very small, was studied by
Forn\ae ss and Sibony in \cite{fsy}. The complement of the forward
Julia set consists of exactly two connected components, called
Fatou components, namely the basin of attraction of the unique
attracting fixed point and the basin of a point at infinity.

It was proved in \cite{wo} that for every $h \in (3,4)$ there
exists a conjugated H\'{e}non map $f$ as above such that the
Hausdorff dimension of the forward Julia set $J^+$ of $f$ is
exactly equal to $h$. In fact, it follows from the proof of
Theorem 4.1 in \cite{wo} that the dimension of $J^+$ is equal to
$h$ in any neighborhood of a point $p \in J^+$.

Now let $h_1, h_2, \ldots$ be an increasing sequence of real
numbers that converge to $4$. Let $f_1, f_2, \ldots$ be a sequence
of H\'{e}non maps as above such that the Hausdorff dimension near
any point of the Julia set of $f_j$ is equal to $\frac{h_j+4}{2}$.

We conjugate the maps $f_j$ with translations that move the
attracting fixed point to zero to get that $f_j$ has an attracting
fixed point at zero. Of course, the Hausdorff-dimension properties
of the maps are unchanged by this conjugation. For every map there
is a small ball centered at 0, say $B(r_j)$, such that $f_j$
attracts to $0$ on $B(r_j)$.

We will denote $J^+_j$ for the forward Julia set of $f_j$. We
choose an increasing sequence $R_1, R_2, \ldots$ such that the
constant $R_j$ defines a filtration $D_j, V^+_j, V^-_j$ for the
map $f_j$, and such that $\lim R_j = \infty$.

In the theorem below we will construct integers $n_1, n_2, \ldots$
and study the forward Julia set of the sequence $f_1^{n_1},
f_1^{n_2}, \ldots$, which we will denote by $J^+$. As before we
will write $F_j$ for $f_j^{n_j}$

\begin{theorem}\label{box}
We can choose integers $n_1, n_2, \ldots$ large enough such that
the upper box dimension of $J^+$ is equal to $4$ near any point in
$J^+$.
\end{theorem}

By large enough we mean exactly the same as in Theorem
\ref{repeat}, namely that we may increase the value of any $n_j$,
but that we then might have to increase the values of all
subsequent integers.

\begin{proof}
First of all we make sure that all $n_j$'s are so large that
$F_j(B(r_j)) \subset B(r_{j+1})$ and such that $F_j(V_j^+) \subset
V_{j+1}^+$. We also make the $n_1, n_2, \ldots$ so large that if
an orbit is such that $z_j \in B(r_{j+1})$ (or in $V_{j+1}^+$)
then the orbit converges to the origin (or to the attracting point
at infinity respectively).

We will now define the sequence $n_1, n_2, \ldots$ inductively.
Suppose that we have fixed $n_1, \ldots , n_{j-1}$. We will
proceed to fix $n_j$.

Let $I_j$ be the forward Julia set of the sequence $F_1, F_2,
\ldots F_{j-1}, f_j, f_j, \ldots$. We have that $z \in I_j$ if and
only if $F(j-1)(x) \in J_j^+$, and thus $I_j = F(j-1)^{-1}J_j^+$.
Therefore we have that the Hausdorff dimension of $I_j$ near any
point of $I_j$ is also equal to $\frac{h_j+4}{2}$.

Let $\mathcal{B}_j$ be a finite collection of open balls of radius
$(\frac{1}{2})^j$ each of which intersects $I_j$ and such that
$\mathcal{B}_j$ covers $I_j \cap D_j$. Since the Hausdorff
dimension of $I_j$ is strictly larger than $h_j$ near any point in
$I_j$, we can choose some small $\hat{\epsilon}_j > 0$ such that
$\gamma_{h_j}^{\hat{\epsilon}_j}(B\cap I_j) > 2^{j+1}$ for any $B
\in \mathcal{B}_j$.

Now let $\epsilon_j < \hat{\epsilon}_j$ such that
$$
\left(\frac{\hat{\epsilon}_j}{\epsilon_j}\right)^{h_j} < \; 2,
$$
and let $\delta_j < \hat{\epsilon}_j - \epsilon_j$.

Define the set $K_j$ by
$$
K_j = \{ z \in D_j \mid d(z, I_j) \ge \delta \}.
$$

For every $z \in K_j$ we have that $f_j^nF(j-1)(z)$ converges
either to $0$ or to the attracting point at infinity, and we can
also see that $K_j$ is compact. Therefore we may choose $n_j \in
\NN$ such that $F(n)(z) = f_j^{n_j}F(j-1)(z)$ lies in $B(r_{j+1})$
or $V^+_{j+1}$ for any $z \in K_j$. It follows from the
compactness of $I_j \cap D_j$ that we can further increase $n_j$
if necessary such that for any $z \in I_j \cap D_j$ there exist
$x, y \in B(z, \delta_j)$ such that $F(n)(x) \in B(r_{j+1})$ and
$F(n)(y) \in V^+_{j+1}$. It follows that with this choice of $n_j$
we have that $d_H(I_j \cap D_j, J^+ \cap D_j) < \delta_j$.

We continue to construct $n_{j+1}, \ldots$ in the same manner.

Now let $z \in J^+$, let $h \in (0,4)$, and let $\epsilon>0$. It
suffices to show that the upper-box dimension of $J^+ \cap B(z,
\epsilon)$ is larger or equal to $h$.

We have that for $j\in \NN$ large enough the following hold: $h_j
> h$, $z \in D_j$ and $\epsilon > 3(\frac{1}{2})^{j}$. Then there
exist a $B \in \mathcal{B}_j$ such that $B \subset B(z,
\epsilon)$. Let $\{B_i\}$ be an $\epsilon_j$ covering of $J^+ \cap
B(z, \epsilon)$. We write $\tilde{B}_i$ for the ball with the same
center as $B_i$ but with radius $\hat{\epsilon}_j$. Since
$\delta_j < \hat\epsilon_j - \epsilon_j$ and $d_H(I_j \cap D_j,
J^+ \cap D_j) < \delta_j$ we have that $\{\tilde{B}_i\}$ is an
$\hat{\epsilon}_j$-covering of $B$, and thus that
$\hat{\epsilon}_j^{h_j} \# \tilde{B}_i > 2^{j+1}$. Since
$\hat{\epsilon}_j^{h_j} / \epsilon_j^{h_j} < 2$ we have that
$$
\gamma_{h}^{{\epsilon}_j}(J^+ \cap B(z, \epsilon)) >
\gamma_{h_j}^{{\epsilon}_j}(J^+ \cap B(z, \epsilon)) > 2^j.
$$

This holds for every $j$ large enough, therefore we have that
$\mu_h(J^+ \cap B(z, \epsilon)) = \infty$, which completes the
proof.
\end{proof}

We can now prove:

\begin{corollary}
There exists a Fatou-Bieberbach domain $\Omega$ whose boundary has
upper box dimension $4$ near any $z \in \Omega$.
\end{corollary}
\begin{proof}
It follows from Theorem \ref{repeat} that we can enlarge the
sequence $n_1, n_2, \ldots$ if necessary to get that the basin of
attraction of the origin in the proof of Theorem \ref{box} is a
Fatou Bieberbach domain. We will call this basin $\Omega$.

In the proof of Theorem \ref{box} we made sure that every open
ball in $\CC^k$ contains a point that either converges to the
origin or to the attracting point at infinity. Therefore there are
only two Fatou components, and one of them is $\Omega$. Therefore
we have that $\partial \Omega$ is exactly equal to $J^+$ and we
are done.
\end{proof}

It would be interesting to know whether the Hausdorff dimension of
the boundary of $\Omega$ is also equal to $4$.

\section{Hausdorff Dimension $4$}

In our construction we will make use of the following version of
Theorem 2.3 in [FR]:
\begin{theorem}\label{FORST}
Let $K_1, K_2,...,K_m$  be pairwise disjoint polynomially convex
compact sets in $\CC^k$  so that their union is polynomially
convex, and assume that $K_1,K_2,..,K_l$ are star-shaped ($l\leq
m$).  Let $\phi_i\in Aut(\CC^k)$ be automorphisms for $i=1,...m$,
so that the sets $K_i'=\phi_i(K_i)$ and the sets $K_{l+1},...,K_m$
are pairwise disjoint, and so that their union is polynomially
convex. Let $\e>0$.  Then there exists an automorphism $\phi\in
Aut(\CC^n)$,  so that $\|\phi(z)-\phi_i(z)\|<\e$  for all $z\in
K_i$  for $i=1,..,l$,  and so that $\|\phi(z)-z\|<\e$  for all
$z\in K_i$  for $i=l+1,..,m$.
\end{theorem}

The proof of this theorem is a small modification of the proof of
Theorem 2.3 in \cite{fr}. We will give a short outline of how one
makes this modification:

Define $K_0=\cup_{j=l+1}^m K_j$.  Choose $B_R$  so that
$K_0\subset B_R$, and let $p_i$ denote the contracting center of
$K_i$ for $i=1,...,l$. By choosing appropriate $\mathcal C^2$
paths from the $p_i$'s to separate points outside of $\overline
B_R$, and by using these to define an isotopy of biholomorphism,
it is clear that by using the arguments from the proof of Theorem
2.3 in [FR], one can find an automorphism that stays close to the
identity in $K_0$, and maps the other $K_i$'s outside of $B_R$.
Now Theorem 2.3 applies, and we can find an automorphism that
moves $K_0$ far away, while staying close to the identity on the
images of the other $K_i$'s. If we have moved $K_0$  far enough,
we may now use the inverses of our paths to map the images of the
other $K_i$'s approximately back to where we started.  We are now
in a situation so that Theorem 2.3 applies to give us a single
automorphism that approximates each $\phi_i$  good on each $K_i$
for $i=1,...m$, and that stays close to the identity on the image
of $K_0$.  Now we repeat the procedure to map the image $K_0$
approximately back to $K_0$, i.e. define paths to move the images
of the $K_i$'s so far away for $i=1,...m$, that we can move the
image of $K_0$ approximately back to $K_0$, and then use the
inverse paths to finish the construction.

\begin{remark}\label{convex}
Recall also these basic facts concerning polynomially convex sets:
(i) The union of a polynomially convex set and a finite set of
points is polynomially convex.  (ii) If $K_1\cup K_2$  is
polynomially convex, if $K_1\cap K_2=\emptyset$, and if
$K_1'\subset K_1$ is polynomially convex, then $K_1'\cup K_2$  is
polynomially convex. (iii) A polynomially convex set has a
neighborhood basis consisting of sets whose closure is
polynomially convex.

One can prove these facts as follows:

To show (i), let $q\in\CC^k\setminus(K\cup\{p\})$, and let
$f\in\mathcal O(\CC^k)$ so that $f(q)=0, f(p)\neq 0$. Since $K$ is
polynomially convex there exists a $g\in\mathcal O(\CC^k)$ so that
$\|g\|_K<1$  and so that $g(p)=1$.  If $m$ is large enough,
$g^m\cdot f$  separates $q$ from $\{p\}\cup K$, and the result
follows by induction.  To prove (ii), let
$p\in\CC^k\setminus(K_1'\cup K_2)$. There exists a $f\in\mathcal
O(\CC^k)$  separating $p$ from $K_1'$.  We get (ii) by defining a
constant function on $K_2$, and applying the Oka-Weil Theorem [RA,
p.220]. A polynomially convex compact set has a neighborhood basis
$U_1\supset\supset U_2\supset\supset ...$ consisting of analytic
polyhedra defined by entire functions, i.e. a Runge and Stein
neighborhood basis [RA, p.71].  This means that the sets
$(\widehat{\overline U_i})^\circ$  satisfy (iii).
\end{remark}

Theorem \ref{FORST}  can be used with Theorem \ref{square} to
prove a slight generalization of Theorem 8.5 in [RR].  This
theorem states that for any sequence of points
$\{p_j\}\subset\CC^2$ and any convex compact set $K$, one can find
a Fatou-Bieberbach Domain $\Omega\subset\CC^k\setminus K$  so that
$\{p_j\}\subset\Omega$. We can exchange polynomially convex for
convex.

\begin{theorem}\label{dense}
Let $K$  be a polynomially convex compact subset of $\CC^2$, and
let $\{p_j\}_{j\in\NN}\subset\CC^2\setminus K$. Then there exists
a Fatou-Bieberbach domain $\Omega$  so that
$\{p_j\}\subset\Omega\subset\CC^2\setminus K$.
\end{theorem}
\begin{proof}
We may assume that $K$  does not intersect the unit ball in
$\CC^2$.  Let A be the linear map defined by
$$
A:(z,w)\rightarrow (\frac{z}{2},\frac{w}{2})
$$
By Theorem \ref{square}, there is a $\d>0$ so that if
$\{f_j\}_{j\in\NN}\subset Aut_0(\CC^2)$ is a sequence of
automorphisms so that $\|f_i-A\|_{\overline B}<\d$  for all
$i\in\NN$, then $\Omega_{\{f_j\}}$  is biholomorphic to $\CC^2$.

We will construct a sequence of automorphisms by induction, and
the following will be our induction hypothesis $I_j$:  We have
constructed automorphisms $\{f_1,...,f_j\}$  so that the following
is satisfied:
\begin{align}
\label{a}&\|f_i-A\|_{\overline B}<\d \ \mathrm{for} \ i=1,...j\\
\label{b}&f(j)(p_i)\subset B \ \mathrm{for} \ i=1,..,j\\
\label{c}&f(j)(K)\subset\CC^2\setminus\overline B
\end{align}

$I_1$  is satisfied by assuming $K$  to be far enough away from
the unit ball, by letting $p_1$  be the origin, and by defining
$f_1=A$.  Now, assume that we have $I_j$.  For any $\mu>0$, by
Theorem \ref{FORST}, there exists a $\psi\in Aut_0(\CC^2)$  so
that $\|\psi-A\|_{\overline B}<\mu$, and so that $\|\psi -
id\|_{f(j)(K)}<\mu$.  And there exists a $\phi\in Aut_0(\CC^2)$ so
that $\|\phi - id\|_{f(j)(K)\cup\overline B}<\mu$, and so that
$\phi(f(j)(p_{j+1}))\in\psi^{-1}(B)$.  Then $I_{j+1}$  is ensured
by letting $\mu$  be small enough and defining $f_{j+1}=\psi\phi$.

Now by (\eqref{a}) and Theorem \ref{square},
$\Omega=\Omega_{\{f_j\}}^0$ is biholomorphic to $\CC^2$, and by
(\eqref{b}) and (\eqref{c}), $\Omega$ is a Fatou-Bieberbach domain
satisfying the claims of the theorem.

\end{proof}
\begin{corollary}
There exists a Fatou-Bieberbach domain $\Omega\subset\CC^2$  so
that the 4-dimensional Hausdorff measure of $\partial\Omega$ is
non-zero.
\end{corollary}
\begin{proof}
Let $D=(\overline D)^\circ\subset\subset\CC$  be a simply
connected open set so that $\partial D$ has non-zero 2-dimensional
Hausdorff measure. Then $K=\overline{D\times D}\subset\CC^2$  is a
polynomially convex compact set whose boundary has non-zero
4-dimensional Hausdorff measure.  Let $P=\{p_j\}$  be dense in
$\CC^2\setminus K$.  By theorem \ref{dense}, there is a
Fatou-Bieberbach domain $\Omega$ so that
$P\subset\Omega\subset\CC^2\setminus K$, which means that
$\partial K\subset\partial\Omega$.
\end{proof}

Now, this corollary tells us only that a part of boundary of
$\Omega$ is large; it does not give us any information about the
size of boundary near other points.  We now prove Theorem
\ref{hausdorff}:
\\
\\
{\bf{Theorem \ref{hausdorff}}} {\emph{ There exists a
Fatou-Bieberbach domain $\Omega$ so that
$\partial\Omega=\partial\overline\Omega$, and so that the
4-dimensional Hausdorff measure of $\partial\Omega$  near any
point is non-zero.}}
\begin{proof}
Let $D=(\overline D)^\circ\subset\subset\CC$  be a connected and
simply connected set whose boundary has non-zero 2-dimensional
Hausdorff measure. Then $K=\overline{D\times D}$  is a
polynomially convex compact set whose boundary has non-zero
4-dimensional Hausdorff Measure. Let $p\in\CC^2$.  We will let
$K_\e(p)$ denote any such $K$  so that $p\in K$  and so that
$K\subset B(p,\e)$. \

Let $A$ be the linear map defined by
$$
A\colon(z_1,z_2)\rightarrow(\frac{z_1}{2},\frac{z_2}{2})
$$
It follows from Theorem \ref{square} that we can choose a $\d>0$
so that if we construct a sequence of automorphisms
$\{f_j\}\subset Aut_0(\CC^2)$  so that $\|f_i-A\|_B<\d$  for all
$i\in\NN$, then $\Omega_{\{f_j\}}$ is a Fatou-Bieberbach domain.
Also, choose a sequence of strictly positive numbers
$\{\e_j\}_{j\in\NN}$ converging to zero.

We are going to construct a sequence of automorphisms inductively,
and the following will be our induction hypothesis $I_j$: We have
automorphisms $\{f_1,...,f_j\}\subset Aut_0(\CC^2)$, and we have a
polynomially convex compact set $K^j=K^j_1\cup ...\cup
K^j_m\subset B(j+1)$, where each $K^j_i$  is a set $K(p,\e)$  for
some point $p$  and an $\e>0$.  We also have a set of points
$\{t_j\}_{j=1}^l\subset B(j+1)$, and we have a ball $B_{R_j}$
 so that $\overline B_{R_j}\cap\overline B=\emptyset$.  In addition, the following is satisfied:
\begin{align}
\label{1}
&\|f_i-A\|_{\overline B}<\d \ \mathrm{for} \ i=1,...,j\\
\label{2}
&B(j+1)\setminus f(j)^{-1}(\overline B)\neq\emptyset\\
\label{3}
&f(j)(K^j)\cap\overline B=\emptyset\\
\label{4}
&f(j)(t_i)\in B \ \mathrm{for} \ i=1,...,l\\
\label{5}
 &\mathrm{For \ any} \ p\in B(j+1)\setminus (K^j)^\circ \
\mathrm{there \ is \ a}
\ t_i \ \mathrm{so \ that} \ \|t_i-p\|<\e_j\\
\label{6}
&\mathrm{For \ any} \ p\in B(j+1)\setminus f(j)^{-1}(B)
\
\mathrm{there \ is \ a} \ K^j_i \ \mathrm{so \ that} \ d(p,K^j_i)<\e_j\\
\label{7}
&f(j)(K^j)\subset B_{R_j}
\end{align}

We will now show how to construct $f_{j+1}$  so that we get
$I_{j+1}$.  It will be clear from the construction how to
construct $f_1$.

Choose a set of points $\{p_i\}_{i=1}^n\subset
B(j+2)\setminus(K^j\cup f(j)^{-1}(\overline B))$, so that for any
point $p\in B(j+2)\setminus((K^j)^\circ\cup f(j)^{-1}(B)$, there
is a $p_i$  so that $\|p-p_i\|<\e_{j+1}$.  Let $q_i=f(j)(p_i)$ for
$i=1,...,n$.  Then Remark \ref{convex} tells us that $\overline
B\cup f(j)(K^j)\cup\{q_i\}_{i=1}^n$ is a polynomially convex
compact set, and that we may choose a $\r>0$ so that the following
is satisfied:
\begin{align}
&\overline B\cup f(j)(K^j)\cup(\cup_{i=1}^n\overline{B(\r)(q_i)})
\ \mathrm{is \ polynomially \
convex}\\
&\overline{B(\r)(q_i)}\cap(\overline B\cup
f(j)(K^j))=\emptyset \ \mathrm{for} \ i=1,...,n\\
&\overline{B(\r)(q_i)}\cap\overline{B(\r)(q_l)}=\emptyset \
\mathrm{for} \ i\neq l.
\end{align}

For each $i$, let $\tilde K_i=K_\r(p_i)$. Define $K^{j+1}=K^j\cup
f(j)^{-1}(\cup_{i=1}^n\tilde K_i)$, which by Remark \ref{convex}
also is polynomially convex. This takes care of (\ref{6}).

Let $\{t_j\}_{j=1}^k\subset B(j+2)\setminus K^{j+1}$  so that for
any $p\in B(j+2)\setminus(K^{j+1})^\circ$, there is a $t_i$  so
that $\|p-t_i\|<\e_{j+1}$.  This ensures (\ref{5}).  Let $\tilde
t_i=F(j)(t_i)$  for $i=1,...,k$. \

We will now construct $f_{j+1}$, and for each statement made on
the existence of a certain automorphism, we refer to Theorem
\ref{FORST}. For any $\mu
>0$ there exists a $\varphi\in Aut_0(\CC^2)$  so that
$\|\varphi(z)-z\|<\mu$ for all $z\in\overline B\cup F(j)(K^j)$,
and so that there exists an $R_{j+1}$ so that $\overline
B_{R_{j+1}}\cap\overline B=\emptyset$, so that $\tilde
K^{j+1}=\varphi(\cup_{j=1}^n\tilde K_i\cup F(j)(K^j))\subset
B_{R_{j+1}}$. There exists a $\phi\in Aut_0(\CC^2)$  so that
$\|\phi(z)-z\|<\mu$ for all $z\in\overline B_{R_{j+1}}$, and so
that $\|\phi-A\|_{\overline B}<\mu$.  Last, there exists a
$\psi\in Aut_0(\CC^2)$  so that $\|\psi(z)-z\|<\mu$  for all
$z\in\overline B\cup\tilde K^{j+1}$, and so that
$\psi(\varphi(\tilde t_i))\in\phi^{-1}(B)$  for $i=1,...k$.  If we
choose $\mu$  small enough and define $f_{j+1}=\phi\psi\varphi$,
we have now ensured (\ref{1}),(\ref{2}),(\ref{3}),(\ref{4}) and
(\ref{7}), so this completes the induction step.

We have constructed a sequence of automorphisms
$\{f_j\}_{j\in\NN}$, and by Theorem \ref{square},  the basin of
attraction to zero is biholomorphic to $\CC^2$.  And since it is
obviously not the whole of $\CC^2$, it is a Fatou-Bieberbach
domain.  We denote it $\Omega$. It is clear from (\ref{3}) that
none of the sets $K^j$ in the above construction is contained in
$\Omega$, and it is clear from (\ref{4}) that all of the $t_i$'s
chosen at every step are. Let $K_i$ be one of the sets at step
$j$, and let $p\in\partial K_i$. Because of (\ref{5}), there is a
sequence of $t_i$'s converging to $p$, and as all of these are in
$\Omega$, we have that $p\in\partial\Omega$.  So the increasing
union $K=\cup_{j=1}^\infty\partial K^j$  is a set whose
4-dimensional Hausdorff measure is non-zero at any point, and we
have $K\subset\partial\Omega$.  It follows from (\ref{6})  that
$K$ is dense in $\partial\Omega$, and this completes the proof.
\end{proof}

\section{Hausdorff dimension 3}

Let $\Omega$ be a Fatou-Bieberbach domain in $\CC^k$. If the
complement of $\Omega$ has non-empty interior then it is easy to
see that the Hausdorff dimension of $\partial\Omega$ is at least
$2k-1$. It is however possible for a Fatou-Bieberbach domain to be
dense in $\CC^k$ (see \cite{rr} or Theorem \ref{dense}), so this
does not guarantee that the dimension of $\partial \Omega$ is at
least $2k-1$. However, we can prove the following:
\begin{theorem}\label{three}
Let $\Omega$ be a Fatou-Bieberbach domain in $\CC^k$ which is
Runge. Then the Hausdorff dimension of $\partial\Omega$ is at
least $2k-1$ near any point of the boundary.
\end{theorem}
\begin{proof}
Assume that $0 \in \partial\Omega$, and we may assume that
$\Omega$ is dense. We will show that the dimension of
$\partial\Omega$ is at least $2k-1$ in any neighborhood of $0$.

For $z \in \CC^k$ we will write $z = (z_1, \ldots ,z_k)$. Let $U$
be an open subset of a small ball centered at the origin
consisting of annuli that are uniformly bounded away from the
hyperplane $z_k = 0$. Define $f : U \rightarrow \CC^{k-1} \times
\RR$ by
$$
f(z) = (\frac{z_1}{z_k} , \ldots, \frac{z_{k-1}}{z_k}, \|z\|).
$$

This is a smooth open mapping, and the preimages of points are
circles centered at the origin. Since $\Omega$ is Runge we have
that every circle centered at the origin intersects
$\partial\Omega$. Therefore we have that $f(\partial\Omega \cap U)
= f(U),$ which is an open subset of $\RR^{2k-1}$. Therefore the
dimension of $U \cup
\partial\Omega$ is at least $2k-1$.
\end{proof}
It is not known whether there exist Fatou-Bieberbach domains that
are not Runge. Therefore Theorem \ref{three} does not rule out the
existence of a Fatou-Bieberbach domain whose boundary has
Hausdorff dimension strictly less than $3$.

\end{document}